\documentclass{article}
\usepackage{amsmath, amsfonts}
\usepackage{amsthm}
\usepackage{color}
\usepackage{graphicx}
\usepackage{algorithm}
\usepackage{algpseudocode}
\usepackage{array}
\usepackage{multirow}
\usepackage{tabularx}
\usepackage{subcaption}
\usepackage{authblk}
\usepackage{booktabs}
\usepackage[table]{xcolor}
\usepackage{array} 
\usepackage{geometry}
\usepackage{rotating}
\usepackage{caption}

\numberwithin{figure}{section}
\numberwithin{table}{section}
\numberwithin{equation}{section}

\theoremstyle{remark}

\theoremstyle{definition}

\theoremstyle{assumption}

\begin{document}

%
%\title{A Methodology for Robust Optimization involving Expensive Function Evaluations with Application to Airfoil Design Optimization under Uncertainty}
\title{Robust Airfoil Design Optimization via a Bilevel Model-Based Methodology}
\author[1]{Alpaslan Kurt}
\author[1]{Dilara Bük}
\author[1]{Yusuf Yazıcıoğlu}
\author[1]{Figen Öztoprak\footnote{Corresponding author.  E-mail: figen.topkaya@gtu.edu.tr}}
\author[2]{Onur Son}
\author[3]{Görkem Demir}
\affil[1]{Gebze Technical University, Industrial Engineering Department, Kocaeli, Turkey}
\affil[2]{Gebze Technical University, Aeronautical Engineering Department, Kocaeli, Turkey}
\affil[3]{Turkish Aerospace (TUSA\c{S}), Ankara, Turkey}
\date{}
\setcounter{Maxaffil}{0}
\renewcommand\Affilfont{\itshape\small}
\maketitle            

\begin{abstract}
We propose GLORO, a bilevel solution methodology for robust optimization involving expensive function evaluations.  The methodology is carefully designed to achieve satisfactory optimization results in a computationally efficient manner.  It is based on Gaussian Process surrogates in the parameter space for approximate computations of the objective / constraints of the robust optimization problem via Bayesian optimization (lower level), and local models constructed in the variable space using these approximate evaluations (upper level).  Both the use of Bayesian optimization (in parameter space) and the use of local models (in variable space) are motivated by the idea of guiding the expensive function evaluations to the regions of interest for the optimization process.  %The methodological work is motivated by and tested on a robust airfoil design optimization problem.  The results obtained for the RAE2822 airfoil are encouraging.
The methodological work is motivated by and tested on a robust airfoil design optimization problem. This application focuses on the RAE 2822 airfoil, optimizing its shape to ensure a robust lift-to-drag ratio under operational uncertainties in Mach number and angle of attack.\\

\smallskip

\end{abstract}

%\keywords{Robust optimization; Gaussian Process; airfoil design optimization.}

\section{Introduction}

Robust optimization has evolved since 1990s as an approach to deal with uncertainties.  The basic idea is to produce solutions less sensitive to the changes in problem parameters by incorporating those parameter uncertainties in the formulation of the optimization problem\cite{mulvey1995robust, ben2002robust,gorissen2015practical}.  A solution to an optimization problem is considered to be robust if it stays close to optimal for any realization of the uncertain problem parameters.  

The work in this paper is motivated by the design optimization problem of an airfoil under uncertainties in operating conditions.  Conventional techniques of robust optimization cannot be directly used for solving this problem because the quantities of interest (the objective and constraint functions of the design optimization problem) can only be computed via computational fluid dynamics (CFD) software; thus, the robust reformulations of these quantities do not have closed analytical forms.  Moreover, the CFD simulations needed for function evaluations are quite computationally expensive which introduces an important restriction for the solution process.  

A popular approach for deterministic airfoil design optimization is to use CFD simulation outputs directly %(without first fitting a global surrogate) 
in running derivative-based optimization solvers, by employing the adjoint computations provided by the CFD software.  However, the usefulness of adjoints for robust formulations is not clear since the approximate evaluations for such formulations can introduce significant magnitudes of noise.  %On the other hand, derivative-free optimization(DFO) literature has limited amount of work on the behavior of methods in the presence of noise, especially for constrained problems.  
There are several studies in the literature that investigate robust modelling, uncertainty quantification, and solution techniques for the aerodynamic shape optimization problem in the presence of uncertainties; e.g. \cite{cook2017robust, papadimitriou2016aerodynamic,fusi2020assessment,tao2019application, nemati2020robust}.  In Appendix B, we present a summary table covering representative work. %,which are mostly based on global surrogates and metaheuristic techniques.  

Generally, robust airfoil optimization requires addressing uncertainty representation, robustness/reliability formulation (whether as objectives or constraints), and computational tractability. Cook and Jarrett\cite{cook2017robust} showed that representing uncertainty as probabilistic or interval based can change which design qualifies as robust.  Papadimitriou and Papadimitriou\cite{papadimitriou2016aerodynamic} minimized a weighted sum of the mean and standard deviation of drag an airfoil under a reliability constraint on lift.  Fusi et al.\cite{fusi2020assessment} compared deterministic and robust formulations for airfoils, finding better off-design performance at the cost of harder trim-constraint handling.  Beyond the common mean-variance trade-off, Croicu et al.\cite{croicu2012robust} compared maximum expected value and expected maximum value criteria, favoring the latter for its ease of implementation with existing solvers.

To evaluate the mean and variance utilized in these formulations without relying on costly Monte Carlo sampling, the polynomial chaos method has been adopted as an efficient uncertainty propagation technique\cite{dodson2009robust}.  More recent variants continue to extend this approach to adjoint-coupled and transonic settings\cite{ji2025d}.  Surrogates other than polynomial chaos are also used for uncertainty propagation in order to model the effect of uncertainties in problem parameters on the resulting function evaluations; i.e. simulation outputs\cite{keane2020robust}.  
%Nevertheless, the repeated CFD evaluations required for polynomial chaos-based uncertainty quantification throughout the optimization process still pose a major computational burden. This has motivated the widespread adoption of surrogate modeling techniques. 
Tao and Sun\cite{tao2019application} used a multi-fidelity deep-learning surrogate to optimize airfoil shapes under Mach number uncertainty.  Keane\cite{keane2012cokriging} proposed cokriging as a general surrogate strategy for robust design, combining correlated low- and high-fidelity data to reduce expensive evaluations.  More recently, Zhang et al.\cite{zhang2024efficient} combined a multi-level hierarchical Kriging model with a multi-fidelity expected improvement criterion, using cheap low-fidelity data to guide where high-fidelity simulations are spent for both optimization and uncertainty quantification.

For both optimization and uncertainty quantification purposes, it is not straigtforward to build successful surrogate models efficiently.  Achieving efficiency in surrogate building requires a careful selection of the sampling points as well as careful decisions regarding the model complexity\cite{cozad2014learning, catalani2023comparative}. 

%An appealing approach is to employ surrogate models to replace expensive function evaluations %as having deterministic differentiable global surrogates
%required during the optimization process.  %However, building successful global surrogate models is not straightforward, and may require evaluations far from the path of optimization\cite{cozad2014learning, catalani2023comparative}. 

In this work, we propose a novel bilevel solution methodology, which we call GLORO, to solve robust optimization problems involving expensive function evaluations.    GLORO relies on surrogate models both for upper-level optimization in the space of design variables, and for lower-level uncertainty propagation taking part in the space of uncertain problem parameters.  To improve sample efficiency, the methodology focuses on the \emph{path of optimization} in building the surrogate models.  For uncertainty quantification, Gaussian Process(GP) surrogates are produced with sampling guided by Bayesian optimization.  On the other hand, local regression models are used for optimization in the space of design variables following the basic steps of model-based derivative-free optimization (DFO).  

There are a few studies in the literature which we believe are relevant to the bilevel solution approach we follow in this work.  The SNOWPAC algorithm by Augustin and Marzouk\cite{augustin2017trust} is a model-based DFO method that uses interpolations.  They use GP to model the noise in robust evaluations, and use it for reducing the approximation error.  Conn and Vicente\cite{conn2012bilevel} propose a two-level model-based DFO method to solve bilevel derivative-free optimization problems.  They consider application of a derivative free method for solving both the upper and lower level problems.  Remarkably, they provide upper bounds on the level of noise that make the existance of noise not interfere with the convergence requirements of the trust region mechanism.  Kolvenbach et al.\cite{kolvenbach2018approach} use DFO for solving an inner optimization problem on the uncertainty set of parameters, which is then employed in an outer level minimization problem in the space of design variables.  Sabater et al\cite{sabater2022robust} use global GP surrogates for both uncertainty propagation and modelling the objective of the design optimization problem.  %However, their surrogate models are global.

The basic novelty of our approach is the integration of Bayesian optimization and model-based DFO techniques in a bilevel scheme.  It builds two levels of models  in the space of uncertain problem parameters and in the space of design variables, respectively.  Both model building phases are embedded into the optimization process to achieve sample efficiency (i.e. in model building, we only focus on the regions of interest for the optimization process).  In particular, we employ Bayesian optimization as a way to efficiently compute approximate robust evaluations of the problem, and use these (noisy) evaluations within an upper level model-based optimization algorithm which relies on randomness to further improve computational efficiency.  %On the other hand, the proposed technique is complementary to the current work in the robust airfoil design optimization literature, which is mostly based on global surrogate building (prior to optimization, not in an embedded way to it) and metaheuristic techniques.  
We elaborate below the distinctions and the contributions of GLORO.
\begin{itemize}
%\item	\textit{Problem formulation.} The problem formulation requires the maximization of the minimum ratio of lift coefficient to the drag coefficient over an uncertainty set of parameters without specifying explicit restrictions on the lift coefficient. [FGN: Bu paragraf tamamlanacak.  Formulasyonla ilgili katkımız hakkında ne soyleyebiliriz?]
\item \textit{Computation of robust evaluations.} %Our uncertainty space is two-dimensional (Mach number and angle of attack are considered as the uncertain parameters).  We employ Bayesian optimization to approximately compute the minimum ratio of lift coefficient to the drag coefficient over that two-dimensional uncertainty set.  
For objective/constraint evaluations of the robust optimization problem, we suggest using Bayesian optimization in the space of uncertain parameters.  The motivation for the use of Bayesian optimization is to provide savings in terms of number of simulations per robust evaluation compared to statistics-based approaches.  As our robust formulation requires minimization in the parameter space, we nicely get the relevant uncertainty quantification output (i.e. the GP model) as a byproduct of Bayesian optimization.  We do not employ a global refinement/infill phase to improve the GP model as the Bayesian optimization process readily guides us to the regions of interest.  We assume that the mapping in between the uncertain parameters and the robust objective/constraints depends on the values of the design variables.  Therefore, we re-run Bayesian optimization each time the design vector is updated (i.e. each time a function evaluation of the robust objective/constraint functions are required for a new design).  The robust evaluations provided by the Bayesian optimization are approximate, which is a fact we take into account in designing the upper level optimization procedure that works in the space of design variables.
\item \textit{Local model based optimization solver.} We follow some of the main ideas in model-based DFO in our upper level optimization technique.  The step computations in the design space are based on local models of the objective and constraint functions of the robust optimization problem.  Using local models avoids the effort of creating a surrogate model in the full solution space.  Instead, it only focuses on the subregions of the solution space around the iterates of the optimization algorithm as these regions are critical for computing successful steps.  In this regard, a local-model based technique needs to produce models of the problem during the optimization process, contrary to the two-stage techniques that produce global surrogate models prior to optimization.  Although GLORO uses local models and trust regions as in model-based DFO methods, it does not fully employ a conventional DFO method as we cannot make as many evaluations as needed to fully satisfy the well-poisedness conditions\cite{conn2009}.  Also unlike the classical model-based DFO methods, the model error does not go to zero when the trust region keeps shrinking.  This is due to the presence of stochastic noise in our robust evaluations resulting from the approximate solutions provided by  Bayesian optimization.  %On the other hand, GLORO is complementary to the existing work on robust design optimization that relies on metaheuristic optimization techniques\cite{arias2011, nemati2020robust, tao2019application}. 
\item\textit{Use of randomness.}  The sample size required by conventional model-based DFO techniques for building local models depends on the problem size.  By sacrificing some accuracy for efficiency, we consider the use of randomness in GLORO to avoid dependence on the problem size.  In the context of model-based DFO, the use of random sample points have been considered in the literature; Bandeira et al.\cite{bandeira2014} provide examples where random sample sets provide better final solutions within less function evaluations compared to deterministic sample sets.  In \cite{maggiar2018} and \cite{chen2015randomized}, randomness is employed via Gaussian smoothing, which is motivated by the existance of noise in function evaluations.  %Remarkably, Maggiar et al. \cite{maggiar2018} employ adaptive multiple importance sampling to reuse all earlier sample points in building local quadratic regression models, and use confidence interval estimations to decide whether new sample points are needed. 
In our case, the robust evaluations contain stochastic noise as mentioned above, but our primary motivation for the use of randomness is scalability.  At each outer iteration of GLORO, we require one additonal evaluation at a random sample point; however, unlike random search methods we use this value together with earlier evaluations to construct a local model and compute an overall step. 
%An alternative for solving the problem at hand would be to use a noise-aware nonlinear algorithm by employing noisy estimates of the derivatives of robust evaluations [Ref].  However, computation of such estimations would unavoidably have dependence on problem dimension.  We therefore propose in this paper a methodology that uses randomness to avoid dependence on problem dimension.
\end{itemize}

In this paper, we utilize our proposed methodology GLORO for solving a robust airfoil design optimization problem.  Our numerical tests verify that it can find acceptable solutions to the robust airfoil design optimization problem within reasonable number of function evaluations (CFD simulations).  The rest of the paper is organized as follows.  In Section~\ref{sec:problem}, we present the mathematical model of the robust optimization problem we study.  Section~\ref{sec:method} contains a detailed description of the proposed methodology, and implementation of this methodology to the robust optimization of RAE2822 airfoil is given in Section~\ref{sec:airfoil}.  Finally, we share our final comments in Section~\ref{sec:future}. 
%The proposed bilevel methodology is a practically promising tool designed carefully for solving a hard optimization problem in an efficient manner.
\section{Problem Statement}\label{sec:problem}

This work is motivated by the 2D design optimization of an airfoil in the presence of uncertainties.  In particular, we aim at computing a robust design that maximizes the lift-to-drag ratio over a set of values of the Mach number and the angle of attack (AoA).  Bound restrictions on the design vector are introduced.  In addition, a constraint requiring a minimum thickness value is introduced similar to the model given in \cite{SU2rae2822}.

Let $x\in \mathbb{R}^n$ denote the design vector, $U_M$ and $U_A$ denote the uncertainty sets for the Mach number and AoA values.  The set of acceptable designs lie within the bounds $[l,u]$; i.e. $x_i\in[l_i,u_i]$ for each design variable $x_i$, $i\in\{1,2\cdots,n\}$.  We define the lift and drag coefficients as functions of the design vector $x$, the Mach number $M$, and the angle of attack $\alpha$, via functions $c_L$ and $c_D$ respectively.  We denote the thickness with $h(x)$, whereas $h_{\min}$ stands for the minimum thickness value allowed.

The resulting robust design optimization problem can be written as follows.
\begin{align}
\label{eq:rdo}
\begin{split}
\max_{x\in \mathbb{R}^n, f\in \mathbb{R}} \ \ & f\\ 
\mbox{such that} \ \  & \min_{M\in U_M, \alpha \in U_A} \frac{c_L(x;M,\alpha)}{c_D(x;M,\alpha)} \geq f\\
& h(x) \geq h_{min}\\
& l\leq x \leq u.
\end{split}
\end{align}
There is indeed a requirement on the pitching moment not being smaller than a threshold.  We do not explicitly impose this as a constraint in our optimization model, but monitor the value of the moment to verify this requirement is readily satisfied. 

Note that \eqref{eq:rdo} is an instance of the following generic robust optimization model
\begin{align}
\label{prob:gro}
\begin{split}
& \max_{x\in \mathbb{R}^n, f \in \mathbb{R}} \ \ f\\ 
\mbox{subject to } \ & c_R(x) \geq f \\
& h(x) \geq 0,  \ \ l\leq x \leq u,
\end{split}
\end{align}
where
\[
c_R(x)=\min_{p\in U_p} c(x;p).
\]
The description of the proposed methodology given in Section~\ref{sec:method} is based on \eqref{prob:gro}, and it has been used specifically for solving \eqref{eq:rdo} in Section~\ref{sec:airfoil}. 

\bigskip

\section{Proposed Methodology}\label{sec:method}

\subsection{Background}

\paragraph{Robust Optimization.}  %Ulbricht makalesi, Bertsimas makalesi.
Robust optimization has been suggested as an alternative approach to deal with uncertainties in problem data.  It aims at computing solutions that are less sensitive to the uncertain parameters of the optimization problem.   Unlike stochastic optimization, no assumption on the distribution of the uncertain parameters is required.  The typical construction of a robust optimization formulation relies on an \emph{uncertainty set}, which contains all possible realizations of the uncertain parameters (i.e. scenarios)\cite{mulvey1995robust}.  Robust formulation of a constraint $c(x;p)\geq 0$ over the uncertainty set of its parameters $U_p$ therefore can be stated as 
\begin{align}\label{eq:robustc}
c(x;p)\geq 0, \mbox{ for all } p\in U_p, 
\end{align}
which contains as many constraints as the cardinality of set $U_p$.  Clearly, the robust formulation leads to a semi-infinite optimization problem when $U_p$ is continuous. Under certain assumptions on the form of the function $c$ and the structure of the set $U_p$, the infinitely many constraints can be reduced to a finite number of constraints\cite{ben2002robust}.  However, such tractable reformulations are not possible for more complicated constraint functions.  Different approaches have been suggested in the literature for this \emph{hard case}; see for instance \cite{zhang2007general, bertsimas2010nonconvex, gorissen2015practical, garreis2017constrained}.  

As for airfoil design optimization, the component functions of the optimization problem are evaluated via CFD simulations, which puts its robust formulation into the hard case category.  On the other hand, robustness measures based on the mean and standard deviation (i.e. a combination of $\mathbb{E}_p[c(x;p)]$ and $\mathbb{V}_p[c(x;p)]$) have often been employed in airfoil design optimization literature rather than the above explained robust reformulation, although both formulations aim at reducing solution sensitivity.  

\paragraph{Uncertainty Quantification.}
For any of the above mentioned robust optimization formulations, the relationship between the uncertain parameters and the robustness measures (i.e. the objective and constraint functions of the robust reformulation) needs to be established.   For robust reformulations based on mean and variance, different ideas have been tried for computing estimations of these statistics; e.g. polynomial chaos\cite{cook2017robust}, Monte-Carlo\cite{ng2014multifidelity}, numerical integration\cite{papadimitriou2016aerodynamic}, GP\cite{dellino2015metamodel}.  

For the robust formulation approach we are interested, it is possible to state the computation of the robustness measure as a lower level optimization problem.  That is, the constraint formulation in \eqref{eq:robustc} can be stated as follows yielding a bilevel optimization problem.
\begin{align}\label{eq:robustc2}
c(x;p)\geq 0, \mbox{ for all } p\in U_p \qquad \equiv \qquad  \min_{p\in U_p} c(x;p)\geq 0.
\end{align}
A remarkable idea is to construct local models for the effect of parameter perturbations on the function values around estimations of the uncertain parameters\cite{zhang2007general}.    Kolvenbach et al.\cite{kolvenbach2018approach} employ local quadratic models of the functions in parameter
space in solving this conventional (worst-case) robust reformulation of the design optimization problem.  However, when the number of uncertain parameters is not large, it is reasonable to tackle \eqref{eq:robustc2} as a global optimization problem, searching for the global minimum value of $c(x;p)$ over all values of $p\in U_p$.  

Note that the local or global models constructed for solving \eqref{eq:robustc2} readily provide maps that model the relationship between $p$ and $c(x;p)$.  When Bayesian optimization is used for solving \eqref{eq:robustc2}, a GP regression model is constructed.  Algorithm~\ref{alg:bo} presents the basic steps of Bayesian optimization when applied for solving the problem
\begin{equation}  \label{eq:probbo}
 \max_{p\in U_p} \ g_x(p) = -c(x;p).
\end{equation} 
A GP model consists of a Gaussian distribution over the function values $g(p)$ with mean function $\mu_g(p)$ and covariance function $\Sigma_g(p)$ on the domain $U_p$.   In Algorithm~\ref{alg:bo}, $D_j$ represents the set of samples accumulated up to iteration $j$.  At each iteration $j$, a sample point $p_j$ is determined with respect to an acquisition function $\alpha$ (e.g. expected improvement), and a new evaluation $g(p_j)$ is computed.  The sample set is then updated by adding the new evaluation, and used to compute the posterior for the GP model.  We refer to \cite{garnett2023} for a detailed discussion of Bayesian optimization.%A GP model consists of a Gaussian distribution over the function values $g(p)$ on the domain $U_p$ with mean function $\mu_g(p)$ and kernel function $\Sigma_g(p)$.  

\begin{algorithm}
\caption{Bayesian optimization framework for solving \eqref{eq:probbo}}
%\caption{Bayesian optimization framework for solving }
\label{alg:bo}
\begin{algorithmic}[1]
\For {$j=1,2,\cdots$}
	\State Update the posterior for GP conditioned on $D_j$ %Update $\mu_g(p)$ and $\Sigma_g(p)$ 
	\State Choose $p_t$ maximizing an acquisition function $\alpha$
\[
p_{j} = \mbox{arg max}_{p\in U_p} \ \alpha(p;D_j)\;
\]  
	\State Compute the function value $g(p_{j})$
	\State $D_{j+1}=D_j\cup (p_j,g(p_{j}))$
\EndFor
\end{algorithmic}
\end{algorithm}

\paragraph{Use of Surrogate Models.} 
We have mentioned above the use of surrogate models (e.g. GP) for uncertainty quantification; i.e. to map uncertain problem parameters to robust objective/constraint values for a given design\cite{sudret2017surrogate}.  In engineering design optimization, surrogate models are also used in the design space to map design variables to objective/constraint values.  Considering the constraint function $c(x;p)$, the former model produces an approximation to $c$ as a function of $p$ whereas the later model produces an approximation to $c$ as a function of $x$.

A surrogate model is generally constructed globally in the full solution space, although the model may be later refined in the regions of interest\cite{forrester2008engineering}.  That typically requires a careful initial sampling in the whole design space.  Another critical issue requiring careful decisions is the complexity of the model class as the form of the underlying function is not known \cite{cozad2014learning}.  On the other hand, surrogate models can be constructed locally during an engineering design optimization process.  This suggests building surrogate models that are good approximations only in regions around the design points obtained as iterates of an optimization algorithm.  Indeed, this is the basic idea of model-based DFO methods.  Use of local modelling allow decreasing the complexity of constructed models and remove the need for an initial sampling on the whole solution space.  However, model construction must be repeated at each iteration of the optimization algorithm.

\paragraph{Presence of Noise.}  As pointed above, the computation of the robust constraint \eqref{eq:robustc2} can be seen as a global optimization problem.  However, global optimization problems are NP-hard -- computation of exact solutions require large computational cost even in small dimensions.  Approximate solutions are reasonable, but they introduce \emph{noise} to the resulting robust constraint evaluation.  More precisely, defining 
\[
c_R(x)=\min_{p\in U_p} c(x;p), 
\]
one can expect computing $c_R(x)$ only approximately.
\[
\tilde c_R(x) \approx c_R(x) \qquad \Rightarrow \qquad \tilde c_R(x) = c_R(x) + \epsilon(x),
\]
where $\epsilon(x)$ is the noise.  When  $\tilde c_R(x)$ is computed via a procedure containing randomness, $\epsilon(x)$ would be random.  %This suggests that our robust optimization problem can be written as a single level optimization problem with noisy evaluations available.
Model-based DFO basically relies on generating and solving local models of the underlying optimization problem, which are conventionally generated via interpolation.  In the presence of noise, one could instead consider using regression models\cite{conn2009} when sufficiently large number of sample points are available.  A recently revoked approach for DFO is to use gradient estimations.  Gradient estimations exploiting randomness have been suggested\cite{nesterov2017random} and have been adapted for the noisy case\cite{chen2015randomized}.  In \cite{chen2018stochastic}, a DFO method based on random models have been proposed to work with evaluations containing stochastic noise.   

Randomness have been exploited also to avoid the bottleneck of dependence on the problem size for model-based DFO; see for instance \cite{menickelly2023avoiding} and \cite{hare2025expected}. 

%For the robust formulation of the problem, we need to model how the uncertainties in flight conditions are propagated to drag and lift computations.  Moreover, we need to model the combined effect of $(m,\alpha)$ onto $c_D$  at given design points $x$ only, ...  As the objective and constraint functions of the robust optimization problem can be computed only approximately, we need to solve an optimization problem with noisy evaluations.  

\bigskip

\subsection{GLORO Algorithm} 
To solve the robust optimization problem in \eqref{eq:rdo}, we propose a bilevel solution methodology.  We consider construction of global surrogates for uncertainty quantification to provide the approximate function evaluations for our robust optimization formulation.  On the other hand, we shift our focus to precision in the design space and suggest the use of local models following ideas from the well-established field of DFO.

Algorithm~\ref{alg:gloro} describes the proposed Gaussian process and LOcal modelling based Robust Optimization (GLORO\footnote{Means \emph{glory} in Esperanto according to Google Translate.}) algorithm for solving the robust optimization model \eqref{prob:gro}.

We let $B(x;\Delta)$ denote the $\ell_\infty$ ball centered at $x$ and has radius $\Delta$.  We also introduce a merit function (for maximization)
\[
\psi(x;\lambda) = c_R(x) +\lambda \min\{h(x),0\},
\]
its local model around an iterate $x_k$
\[
m^\psi_k(x;\lambda,x_k) = m^c_k(x;x_k)+\lambda \min\{m^h_k(x;x_k),0\},
\]
and its available approximate computation (i.e. the noisy evaluation of $\psi$)
\[
\tilde \psi(x;\lambda) = \tilde c_R(x)+\lambda \min\{h(x),0\}.
\]
We do not have a merit parameter update strategy in the algorithm for simplicity.  We assume that a large enough value of $\lambda>0$ is chosen. 

We consider polynomial local models with a $q$-dimensional monomial basis.  In particular, we consider quadratic models in our implementation so that 
\[
q = \frac{1}{2}(n+1)(n+2).
\]
The matrix $\Phi_k$ holds the values of these monomials for the set of data points to be employed in model generation at iteration k, and the vector $y_k$ hold corresponding function evaluations.  In our setting, the number of data points would generally be small compared to $q$, yielding an underdetermined regression problem.  In this case, we construct the minimum $\ell_2$-norm solution of the underdetermined system.  When the number of samples is larger than the dimension of the monomial basis, we suggest using regression models as our approximate robust evaluations contain noise.  We should mention however that for the setting of robust airfoil design optimization, each approximate robust evaluation is so costly that we never expect encountering this later case. 
% When the linear system is underdetermined, then the sklearn.linear_model.LinearRegression finds the minimum L2 norm solution, i.e. argmin_w l2_norm(w) subject to Xw = y
%n-variate 2nd degree polynomial için (n+0-1// 0)+(n+1-1// 1)+(n+2-1// 2) terim var monomial basiste

\begin{algorithm}[htp!] 
\caption{GLORO}
\label{alg:gloro}
\small
\begin{algorithmic}[1]
%\Input 
\Require $x_0$, $U_p$, $\Delta_0$, $N_0\leq q$, $\lambda>0$, $\theta\geq 1$.
%\EndInput

\bigskip

\State \textit{[Initialize.]} Set $k=0$, $X=\emptyset$, $\Phi_0=[]$, $y^c_0=[]$, $y^h_0=[]$.  Generate $N_0$ sample points $x^i$, $i\in\{1,2,\cdots, N\}$ in $B(x_0;\Delta_0)\cap [l,u]$, and rows 
\[
\phi(x^i) := [\phi_0(x^i) \ \ \phi_1(x^i) \ \ \cdots \ \ \phi_q(x^i)]
\]
to $\Phi_0$ such that $\Phi_0$ has full row rank.  $X=X\cup\{x^i\}$, $i\in\{1,2,\cdots, N_0\}$.  Compute approximate robust evaluations 
\[
\tilde c_R(x^i) \approx  \min_{p\in U_p} c(x^i;p)
\]
via Bayesian optimization, for $i\in\{1,2,\cdots, N_0\}$.  Set 
\[
y^c_0 = [\tilde c_R(x^1) \ \ \tilde c_R(x^2) \ \ \cdots \ \ \tilde c_R(x^{N_0})]^T, \ \ \mbox{and} \ \ y^h_0 = [h(x^1) \ \ h(x^2) \ \ \cdots \ \ h(x^{N_0})]^T.
\]

\bigskip

\State \textit{[Model Generation.]}  Generate local models $m^c_k(x;x_k)$ and $m^h_k(x;x_k)$ based on data points $(\Phi_k, y_k^c)$ and $(\Phi_k, y_k^h)$, respectively.

\bigskip

\State \textit{[Step Computation.]}  Solve the local model 
\begin{align}
\label{eq:localmodel}
\begin{split}
x_c  = \ & arg\max_{x\in \mathbb{R}^n} \ \ m^c_k(x; x_k)\\ 
& m^h_k(x) \geq 0,  \ \ l\leq x \leq u,\\
& \|x-x_k\|\leq \Delta_k
\end{split}
\end{align}
to compute a candidate iterate $x_c$.

\bigskip

\State \textit{[Step Assessment.]}  Compute 
\[
\tilde c_R(x^c) \approx \min_{p\in U_p} c(x^c;p)
\]
via Bayesian optimization.  Compare the actual improvement in the merit function with the improvement predicted by its model via computing the ratio
\begin{equation}\label{eq:ratio}
r_k = \frac{\tilde \psi(x^c;\lambda)-\tilde \psi(x_k;\lambda)}{m^\psi_k(x^c;\lambda,x_k)-m^\psi_k(x_k;\lambda,x_k)}
\end{equation}
If $r_k\geq \rho_1$ then $x_{k+1}=x^c$.  Else $x_{k+1}=x_k$. 

\bigskip
 
\State \textit{[Updates.]}  $\Delta_{k+1}$=UpdateTrustRegionRadius($\Delta_k$).

\noindent
$k=k+1$.  $X = X \cup \{x^c\}$.  $\Phi_k=[]$, $y^c_k=[]$, $y^h_k=[]$.
\noindent
For each $x\in X$, if $x\in B(x_k;\theta\Delta_k)\cap [l,u]$ then add row $\phi(x)$ to $\Phi_k$, $\tilde c_R(x)$ to $y^c_k$, and $h(x)$ to $y^h_k$.

\noindent
If $rank(\Phi_k)<q$, generate a new random sample point  $x^r \in B(x_k;\Delta_k)\cap [l,u]$ such that the rank of $\Phi_k$ increases by one.  $X=X\cup\{x^r\}$.  Compute 
\[
\tilde c_R(x^r) \approx \min_{p\in U_p} c(x^r;p)
\]
via Bayesian optimization.  Update $y^c_k$ and $y^h_k$ by adding $\tilde c_R(x^r)$ and $h(x^r)$, respectively.    Go to Step 2.

\end{algorithmic} 
\end{algorithm}

\begin{algorithm}[h!] 
\caption{UpdateTrustRegionRadius}
\label{alg:tru}
\begin{algorithmic}[1]
\Require $\gamma>1$, $0<\rho_1\leq \rho_2<1$, $\Delta_{\min}$.
\If {$r_k\geq \rho_2$}
    \State $\Delta_{k+1} = \gamma \Delta_k$
\ElsIf {$r_k \leq \rho_1$} 
    \State $\Delta_{k+1} = \max\{1/\gamma \Delta_k, \Delta_{\min}\}$
\EndIf 
\end{algorithmic} 
\end{algorithm}

\begin{figure}[htp]
\caption{Evaluation of the constraint value for the robust formulation at a candidate design point.}
\centering
\includegraphics[width=0.6\textwidth]{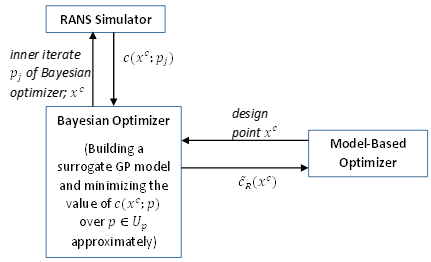}
\label{fig:robusteval}
\end{figure}

The algorithm starts by sampling in $B(x_0;\Delta_0)$.  Initial local models are constructed based on the function evaluations at these initial sample points.  Each time a robust evaluation is needed at a new sample point, Bayesian optimization procedure is called to provide an approximate evaluation.  This process is illustrated in Figure~\ref{fig:robusteval}.  

At each outer iteration,  the local model \eqref{eq:localmodel} is solved within the trust region for computing a step in the space of design variables $x$.  The next iterate is determined with respect to the value of the ratio \eqref{eq:ratio}, which is a measure to assess the improvement in the merit function provided by the candidate design.  The trust region update given in Algorithm~\ref{alg:tru} also depends on that ratio.  We have a lower bound on the trust region radius not to let it converge to zero due to the effect of noise.  The candidate point is always added to the set of sample points.  In addition to the evaluation at this candidate point, one additional evaluation is done at a randomly generated sample point if the model fitting problem is underdetermined in the selected neighborhood of the current iterate.  The additional random sample point is generated within the current trust region, and the random generation is repeated until the rank condition is satisfied.  Then, objective and constraint evaluations are computed at this new sample point.
% Eğer hesaplanan adım sıfıra eşitse lineer modele düşüp onu hesaplıyoruz.

\bigskip

\subsection{Preliminary Tests}

\begin{figure}[htp]
\caption{Value of the constraint $c(x;\mu)$ for different realizations of $\mu_1$ and $\mu_2$.  The cells corresponding to a constraint violation are marked in red.}
\label{fig:toy}
\centering
\includegraphics[width=0.85\textwidth]{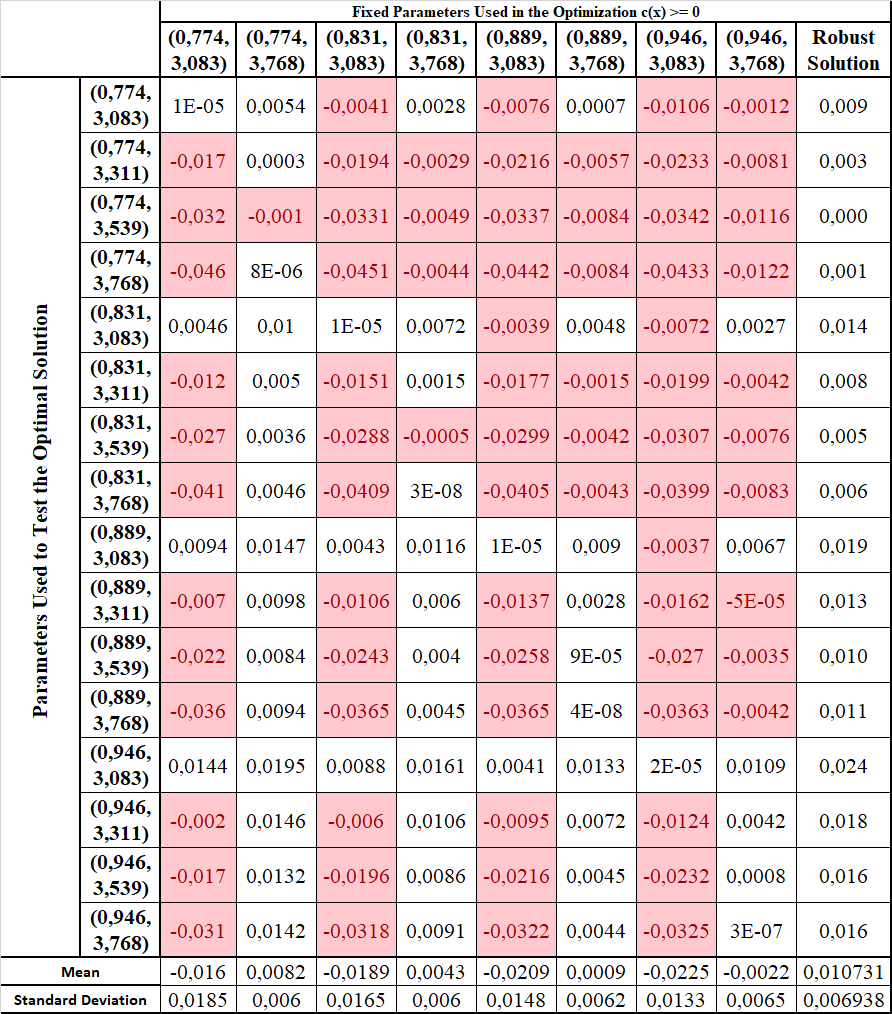}
\end{figure}

We first test GLORO on a toy function with known analytical form.  The problem is the HS91 instance from the Hock \& Schittkowski test set \cite{hock1980} which has a single inequality constraint.  The objective of the problem is a simple convex quadratic function (square of the $\ell_2$-norm of the variables $x$).  However, the constraint function is highly nonlinear as can be seen from its definition given below.
\begin{align*}
c(x;\mu) := &  -\sum_{i=1}^{30}
	\sum_ {j=i+1}^{30} 
	    \mu_i^2 \mu_j^2 A_i A_j \rho_i \rho_j
	    (\sin(\mu_i+\mu_j)/(\mu_i+\mu_j) + \sin(\mu_i-\mu_j)/(\mu_i-\mu_j)) \\
	&+ 0.0001 -\sum_{j=1}^{30} \mu_j^4 A_j^2 \rho_j^2 (\sin(2\mu_j)/(2\mu_j) + 1)/2 \\
	& +\sum_{j=1}^{30} \mu_j^2 A_j \rho_j (
		2\sin(\mu_j)/\mu_j^3 - 2*\cos(\mu_j)/\mu_j^2 ) + 2/15 \geq 0,
\end{align*}
where
\[
\rho_j = -\left(
	    \exp\left(-\mu_j^2 \sum_{i=1}^5 x_i^2\right)
	    +\sum_{t = 2}^n 
		 2 (-1)^{t-1} \exp\left(-\mu_j^2 * \sum_{t=1}^n x_i^2\right)
	    -1\right)/\mu_j^2,
\]
and
\[
A_j = 2\sin(\mu_j)/(\mu_j + \sin(\mu_j)\cos(\mu_j)).
\]
Parameters $\mu_1$ and $\mu_2$ were found to have the most significant impact on the constraint function value $c(x;\mu)$, and therefore were chosen to be the parameters that have uncertainty.  GLORO algorithm is used to solve the problem by setting the uncertainty set to be $\pm 10\%$ of the original values of these parameters; i.e. $\mu_1\in [0.77,0.95]$ and $\mu_2\in [3.08,3.78]$.  Since the constraint has the form 
\[
c(x;\mu)\geq 0,
\]
we expect a robust solution $x$ to provide a nonnegative $c(x;\mu)$ value for all $\mu_1\in [0.77,0.95]$ and $\mu_2\in [3.08,3.78]$, while minimizing the objective function value.  We sampled 100 values of $\mu_1$ and $\mu_2$ from the corresponding uncertainty sets, computed the constraint value provided by the robust solution of GLORO for these 100 parameter realizations, and compared it to the constraint values provided by the optimal solutions of the problem obtained by fixing  the two parameter values to each of the 100 realizations.  We also computed the standard deviation of the $c(x;\mu)$ values for each solution to give an idea on the sensitivity of the solution to these parameters.  A representative subset of the results are given in Figure~\ref{fig:toy}.  The robust solution satisfies the constraint for all realizations of the parameters as desired, whereas optimal solutions obtained via fixed values of these parameters can be infeasible for some realizations.

\bigskip

\section{Robust Optimization of  RAE2822 Airfoil} \label{sec:airfoil}

%We concentrate on the performance of the airfoil within a relatively narrow interval of uncertainty in the hard flight conditions regime.
A two-dimensional RAE2822 airfoil at transonic conditions was analyzed by solving the Reynolds-Averaged Navier–Stokes (RANS) equations coupled with the Spalart–Allmaras (SA) turbulence model.  The open-source computational fluid dynamics (CFD) software SU2\cite{su2} was employed to compute the flow field.  The mesh, configuration, and solver files provided by the SU2 development team were utilized as the basis for the simulations\cite{SU2rae2822}.

\begin{figure}[htp]
\caption{Mesh structure around the RAE 2822 airfoil. (a) full computational domain, (b) near-wall region around the airfoil surface.}
\label{fig:mesh}
\centering
\includegraphics[width=\textwidth]{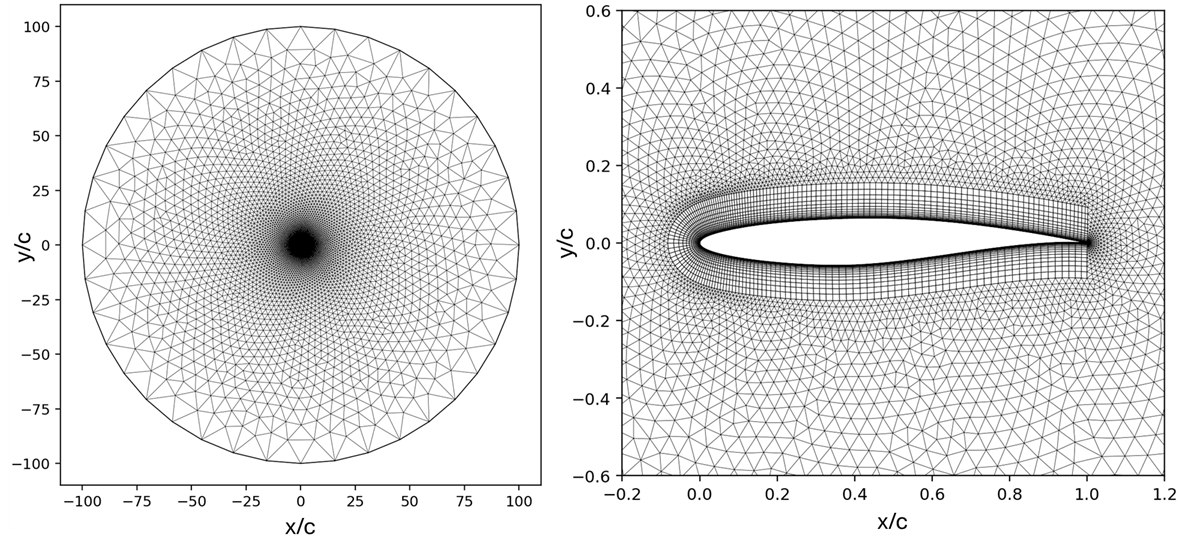}
\end{figure}

The computational mesh is an unstructured O-grid constructed around the RAE2822 airfoil as given in Figure~\ref{fig:mesh}. It consists of 22,842 elements, with the airfoil surface represented by 192 edges and the far-field boundary by 40 edges. The hybrid mesh employs quadrilateral elements in the near-wall region close to the airfoil and triangular elements in the outer domain. The first grid point from the airfoil surface is positioned at a distance of $10^{-5}$ chord lengths, while the far-field boundary is placed 100 chord lengths away.

A wall-resolved approach was employed to maintain a dimensionless wall distance ($y^{+}$) approximately unity, ensuring direct resolution of the viscous sublayer. The airfoil surface was modelled as a no-slip, adiabatic wall. The spatial discretization of the mean flow equations used the Jameson-Schmidt-Turkel (JST) scheme\cite{jameson1981}, while the turbulent terms were discretized using a scalar upwind scheme. Gradients were computed using the weighted-least-squares method, and the Venkatakrishnan slope limiter was applied. For time integration, Euler implicit method was used for both the mean flow and turbulence equations.

\paragraph{Validation.} The \emph{Case 6} configuration from the AGARD Report-138\cite{cook1979aerofoil} is adopted as the validation case for the RAE2822 airfoil. In this case, experiments were conducted at a Mach number of $M=0.725$, a Reynolds number of $Re=6.5\times 10^6$, and an angle of attack of $\alpha=2.92^{\circ}$, yielding a lift coefficient of $C_L=0.743$. The validation simulations were performed at the same Mach and Reynolds numbers, with the angle of attack adjusted to match the target lift coefficient of $C_L=0.743$. This approach has been widely adopted in previous studies to account for wind tunnel wall interference effects \cite{garbaruk2003numerical, da2020sensitivity, zhou1995turbulent, chen2024wind}.  In the present analysis, the target lift coefficient was achieved at an angle of attack of $\alpha=2.44^{\circ}$.

\begin{figure}[htp]
\caption{Pressure coefficient ($C_p$) distribution over the RAE2822 airfoil at $M = 0.725$, $Re = 6.5 \times 10^6$, and $C_L = 0.743$. Experimental results from \emph{Case 6} of Cook et al.\cite{cook1979aerofoil} at $\alpha = 2.92^{\circ}$ (red) are compared with the present numerical results at $\alpha = 2.44^{\circ}$ (blue).}
\label{fig:pr}
\centering
\includegraphics[width=0.8\textwidth]{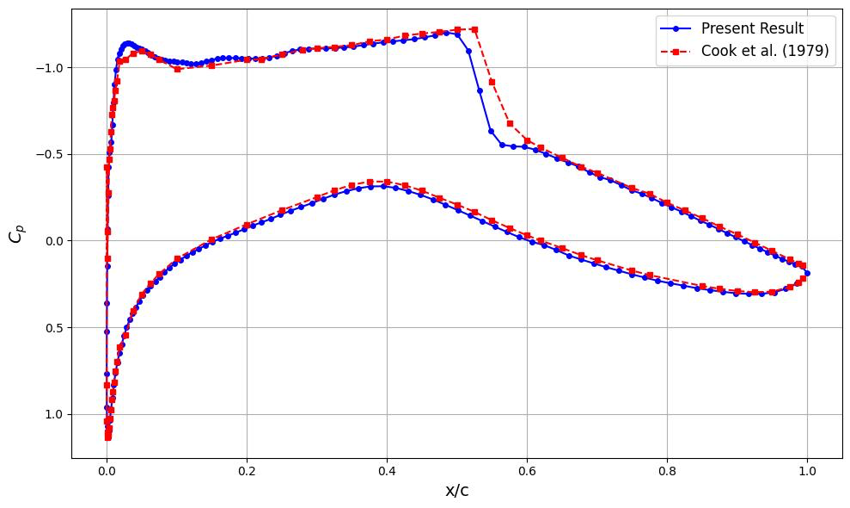}
\end{figure}

The pressure coefficient ($C_p$) distribution over the RAE2822 airfoil is shown in Figure~\ref{fig:pr}. On the upper surface, the distribution exhibits a jump in $C_p$ at approximately $x/c\approx 0.55$, signifying a sudden pressure increase associated with the presence of a shock wave. In the present analysis, the shock position is slightly displaced compared to the reference data, and the suction peak near the leading edge is predicted to be slightly higher. Despite these differences, the computed $C_p$ distribution closely matches the reference experimental data, successfully capturing both the shock-induced pressure jump and general pressure trends along the airfoil surface. 

In addition to the $C_p$ distributions, the corresponding angles of attack required to achieve the same lift coefficient under similar flow conditions are summarized in Table~\ref{tab:cp}. The comparison demonstrates that the present results are in close agreement with previously reported data \cite{cook1979aerofoil, da2020sensitivity, zhang2020inverse, chang1988further}, confirming the validity of the current approach in reproducing both the aerodynamic coefficients and the flow features associated with the RAE2822 airfoil at transonic conditions.

\begin{table}[htp]
\caption{Comparison of present validation results with literature for a fixed lift coefficient of $C_L=0.743$.}
\label{tab:cp}
\small
\centering
\begin{tabular}{m{2.3cm} | c c c | m{1.3cm} m{1.2cm} m{1.3cm}}
\hline
& $C_L$ & $C_D$ & $C_M$ & Angle of attack (deg)	& Mach number	& Reynolds number\\
\hline
Present validation data & $0.743$ & $0.0131$ & $-0.0895$ & $2.44$ & $0.725$ & $6.5 \times 10^6$\\
\hline
Cook et al. (1979) \cite{cook1979aerofoil} & 0.743	& 0.0127 & $-0.095$ & $2.92$ & $0.725$ & $6.5 \times 10^6$ \\
Zhang et al. (2020) \cite{zhang2020inverse} &	0.743	& 0.01219	& - &	2.31	& 0.725 &	$6.5 \times 10^6$\\
Da Ronch et al. (2020) \cite{da2020sensitivity} &	0.743 &	0.0150& 	$-0.0909$	& 2.51	& 0.729	& $6.5 \times 10^6$\\
Chang et al. (1988) \cite{chang1988further} & 0.743	& 0.0112	& - &	2.53& 	0.729	& $6.5 \times 10^6$\\
\hline
\end{tabular}
\end{table}

Following the validation study, attention was directed to the baseline \emph{Case 6} configuration from \cite{cook1979aerofoil}, defined by $M=0.725$ and $\alpha=2.92^\circ$. CFD simulations were conducted for this case prior to initiating the optimization process, and the results were compared with data available in the literature. The summarized findings in Table~\ref{tab:comp} show that the computed lift and drag coefficients are in close agreement with previously reported values at the same Mach number and angle of attack \cite{garbaruk2003numerical, weinmeister2024, vuruskan2019}.  In addition, the pitching moment coefficient was obtained as $C_M=-0.0903$ for this configuration.

\begin{table}[htp]
\caption{Comparison of present results with literature for a fixed angle of attack of $\alpha=2.92^\circ$.}
\label{tab:comp}
\small
\centering
\begin{tabular}{m{2.5cm} | c c c | m{1.3cm} m{1.2cm} m{1.3cm}}
\hline
& $C_L$ & $C_D$ & $C_M$ & Angle of attack (deg)	& Mach number	& Reynolds number\\
\hline
Present RAE2822 airfoil & $0.827$ & $0.0169$ & $-0.0903$ & $2.92$ & $0.725$ & $6.5 \times 10^6$\\
\hline
Garbaruk et al. (2003) \cite{garbaruk2003numerical} & 0.806	& 0.0177 & $-$ & $2.92$ & $0.725$ & $-$ \\
Weinmeister and Sanjaya (2024) \cite{weinmeister2024} & 0.819	& 0.0161& $-$ &	2.92	& 0.725 &	$6.5 \times 10^6$\\
Vuruskan and Hosder (2019) \cite{vuruskan2019} &	0.824	&  0.0209 & $-$ & 2.92	&  0.734	& $6.5 \times 10^6$\\
\hline
\end{tabular}
\end{table}

\begin{figure}[htp]
\caption{Comparison of the original RAE2822 airfoil and the optimized airfoil shapes.}
\label{fig:shpe}
\centering
\includegraphics[width=\textwidth]{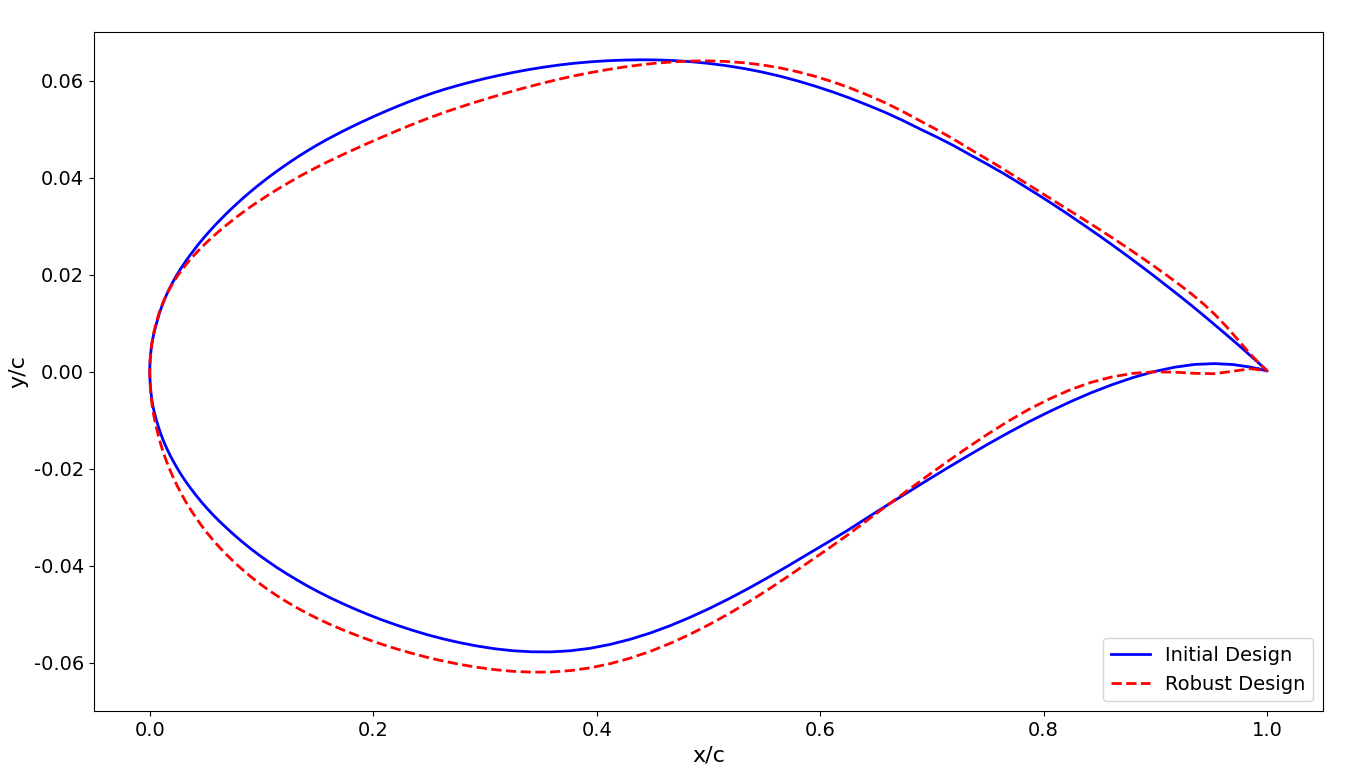}
%\begin{subfigure}{0.85\textwidth}
%\includegraphics[width=\textwidth]{updated_figures/shape}
%\caption{Exact}
%\end{subfigure}
%\hfill
%\begin{subfigure}{0.9\textwidth}
%\includegraphics[width=\textwidth]{updated_figures/bombe_robust}
%\caption{Stretched}
%\end{subfigure}
\end{figure}

\paragraph{Optimization.}  The proposed robust design methodology GLORO is run for robust optimization of the airfoil based on the formulation in \eqref{eq:rdo}.  The design variables are the 38 Hicks-Henne bump function parameters specified in \cite{SU2rae2822} (a chord of 1.0 is assumed for the Hicks-Henne bumps) -- i.e. the size of the problem in the space of design variables is $n=30$.  In Algorithms~\ref{alg:gloro} and \ref{alg:tru}, the initial point is set to $x_0=0_{n\times 1}$, and the parameter values are set to $N_0=10$, $\Delta_0=0.5\times 10^{-5}$, $\lambda=10^4$, $\gamma=2$, and $\rho_1=0.01, \rho_2=0.6$.  Second degree polynomial regression models were fit, and the model degree was reduced to linear as a fallback if no successful steps can be computed via the quadratic model before the trust region value becomes as small as $\Delta_{\min}=10^{-7}$.  To account for variations in the flow conditions, uncertainty ranges were introduced as $\pm 0.01$ for the Mach number and $\pm 0.1^\circ$ for the angle of attack.  Accordingly, the optimization domain was defined as [0.715, 0.735] for the Mach number and [$2.82^{\circ}$, $3.02^{\circ}$] for the angle of attack.  Thus, the Bayesian optimizer employed for computing approximate robust evaluation values solves two-dimensional problems in the parameter space.  Two points were used for initialization of the GP model followed by five Bayesian optimization iterations.  Overall, the optimization process has taken 10 outer iterations of GLORO, and has required a total of 217 SU2 simulations.

The resulting optimized airfoil geometry, obtained within the above mentioned uncertainty bounds, is compared with the original RAE2822 profile in Figure~\ref{fig:shpe}. The figure highlights the geometric modifications introduced by the optimization process, particularly in the suction-side and pressure-side curvatures near the leading-edge region.  Note that the axis scales are adjusted (not fixed) to make these changes more visible.

\begin{sidewaysfigure}[htp]
\caption{Comparison of pressure coefficient distributions between the baseline and optimized airfoils at Mach numbers $M=0.715, 0.725, 0.735$ and for angles of attack $\alpha=2.82^{\circ}, 2.92^{\circ}, 3.02^{\circ}$.}
\label{fig:optp}
\centering
\includegraphics[width=\textwidth]{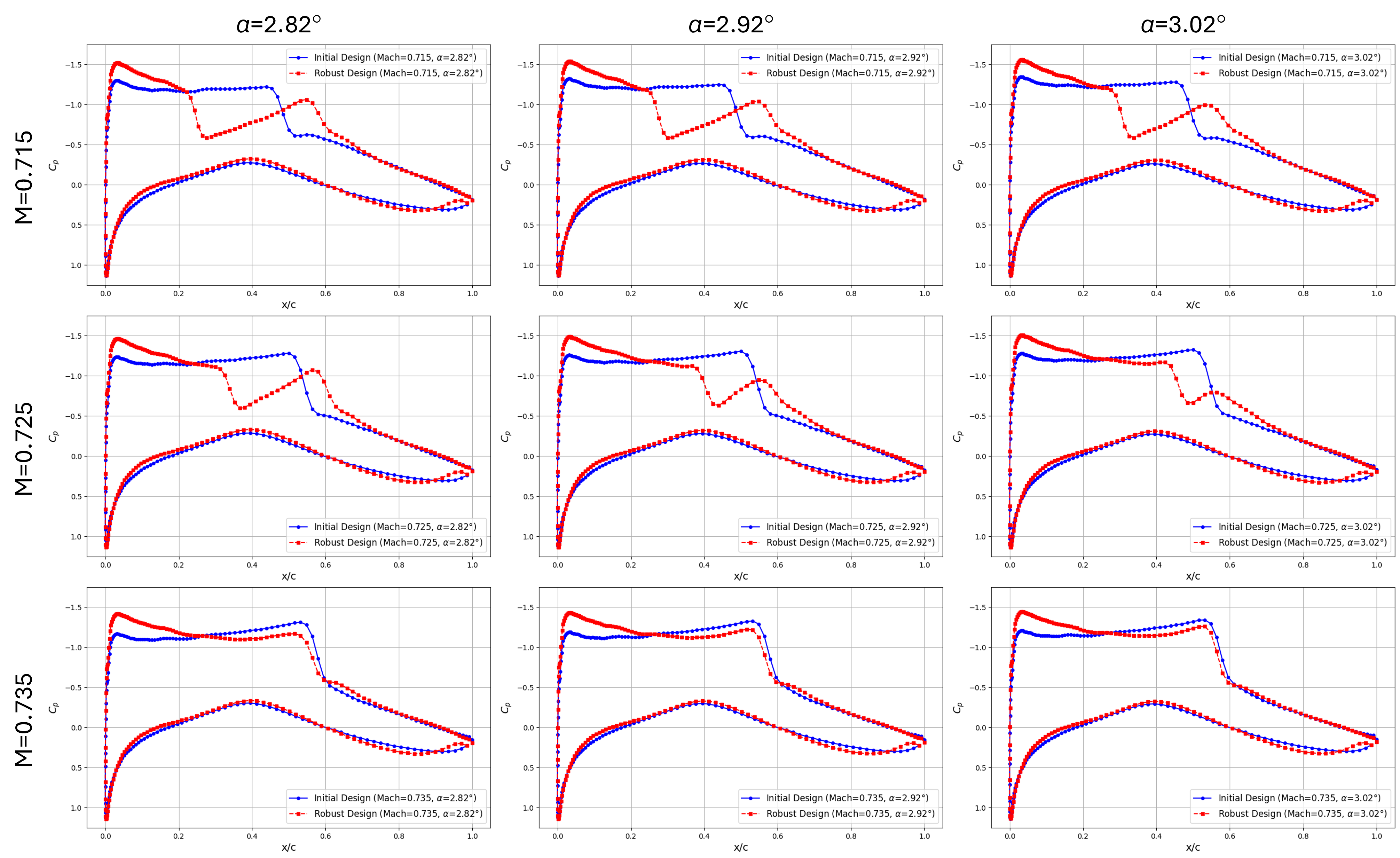}
\end{sidewaysfigure}

The pressure coefficient ($C_p$) distributions for the baseline and optimized airfoil, evaluated at the reference values and within the defined uncertainty bounds, are shown in Figure~\ref{fig:optp}. Compared to the baseline airfoil, the optimized geometry exhibits significant modifications in the surface pressure distribution. At Mach numbers of 0.715 and 0.725, the primary shock wave on the upper surface is shifted noticeably upstream toward the leading edge. Downstream of this shock, the flow undergoes reacceleration, giving rise to a secondary compression region followed by a more gradual pressure recovery.

At the higher Mach number of $M=0.735$, the flow expansion after the first shock and re-compression is no longer evident. Instead, a single, stronger shock forms, with its location closely matching that of the baseline RAE2822 case.  Across all examined conditions, the optimized airfoil generates a higher suction peak near the leading edge in the region $0<x/c<0.20$, which enhances lift production in this region.

The aerodynamic force coefficients corresponding to the pressure distributions in Figure~\ref{fig:optp} are summarized in Table~\ref{tab:compopt}. Cases that exhibit a secondary compression feature in the $C_p$ distribution show a modest reduction in lift relative to the baseline airfoil, where $C_L$ decreasing by approximately 2-6\%. However, in these same cases, drag reduction becomes evident as the Mach number and angle of attack are increased. For example, at $M=0.725$ and $\alpha=3.02^{\circ}$, the drag coefficient decreases by 26\%, while the lift-to-drag ratio improves by 32\% compared to the baseline configuration. The primary shock shifts downstream and secondary compression peak flattens in the cases with drag reduction. 

The emergence of a secondary compression structure in the optimized airfoil resembles to the flow physics associated with shock control bumps (SCBs), which are commonly employed for transonic drag reduction and flow control \cite{zhang2020inverse, jinks2018optimisation,	sabater2020efficient, mazaheri2015}. Interestingly, in the present study, a similar mechanism emerges naturally in the optimized airfoil geometry, without the introduction of a physical bump. This suggests that careful geometric tailoring alone can trigger shockwave interactions analogous to flow control devices, providing drag reduction while maintaining favorable lift characteristics.

For the case of $M=0.735$, where a single shock structure is observed similar to the baseline configuration, a modest lift increase of approximately 0.5\% is obtained. Despite this small change in $C_L$, the drag reduction is substantial, reaching up to 34\% at $\alpha=2.82^{\circ}$.  Consequently, the lift-to-drag ratio improves dramatically, with gains as high as 52\% compared to the baseline airfoil.

%\newgeometry{a4paper,left=1in,right=1in,top=1in,bottom=1in,nohead}
%\begin{landscape}
%\begin{figure}[htp]
%\end{figure}
%\end{landscape}
%\restoregeometry

\begin{table}[htp]
\caption{Comparison of lift coefficient, drag coefficient, and lift-to-drag ratio between baseline and optimized airfoils.}
\label{tab:compopt}
\centering
%\footnotesize
\resizebox{\columnwidth}{!}{
\begin{tabular}{c c |c c c| c c c| c c c}
\hline
& & \multicolumn{3}{c|}{Baseline}& \multicolumn{3}{c|}{Optimized}& \multicolumn{3}{c}{Change(\%)}\\
\hline
$M$ & $\alpha$ & $C_L$ & $C_D$ & $C_L/C_D$ & $C_L$ & $C_D$ & $C_L/C_D$& $\Delta C_L$ & $\Delta C_D$ & $\Delta (C_L/C_D)$\\
\hline
\multirow{3}{2em}{0.715} & 2.820 & 0.8000 & 0.0129 & 62.0691 & 0.7462 & 0.0134 & 55.6826 & -6.72& 3.98 & -10.29\\
& 2.920 & 0.8177& 0.0136 & 60.2150 & 0.7667 & 0.0137 & 56.0678 & -6.24 & 0.70 & -6.89\\
& 3.020 & 0.8354 & 0.0144 & 57.9769 & 0.7884 & 0.0140 & 56.4614 & -5.62 & \textbf{-3.09} & -2.61\\
\hline
\multirow{3}{2em}{0.725} & 2.820 & 0.8104 & 0.0159 & 50.9074 & 0.7674 & 0.0136 & 56.4614 & -5.31& \textbf{-14.62} & \textbf{10.91}\\
& 2.920 & 0.8269& 0.0169 & 48.8339 & 0.7983 & 0.0132 & 60.3319 & -3.46 & \textbf{-21.86} & \textbf{23.52}\\
& 3.020 & 0.8427 & 0.0181 & 46.6622 & 0.8255 & 0.0133 & 61.9235 & -2.04 & \textbf{-26.18} & \textbf{32.71}\\
\hline
\multirow{3}{2em}{0.735} & 2.820 & 0.8087 & 0.0202 & 40.0244 & 0.8133 & 0.0133 & 61.1622 & \textbf{0.56} & \textbf{-34.19} & \textbf{52.81}\\
& 2.920 & 0.8227& 0.0214 & 38.4543 & 0.8260 & 0.0143 & 57.7781 & \textbf{0.41} & \textbf{-33.18} & \textbf{50.25}\\
& 3.020 & 0.8357 & 0.0227 & 36.8905 & 0.8375 & 0.0154 & 54.2501 & \textbf{0.22} & \textbf{-31.85} & \textbf{47.06}\\
\hline
\end{tabular}
}
\end{table}

\begin{figure}[htp]
\caption{Performances of the designs for different parameter realizations.}
\label{fig:var}
\centering
\includegraphics[width=0.6\textwidth]{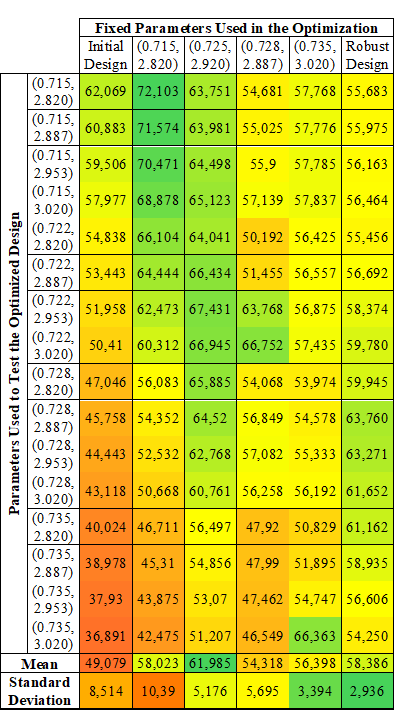}
\end{figure}

It should be emphasized that the objective of this optimization problem is to maximize the lift-to-drag ratio not for single values of the parameters $M$ and $\alpha$, but rather for the entire range of a prescribed uncertainty set.  In Figure~\ref{fig:var}, we contrast the variance of the lift-to-drag ratio for solutions obtained with fixed values of $M$ and $\alpha$ to that of the solution obtained by GLORO.  The solutions with fixed parameter values are obtained via the same outer model-based technique as in GLORO; however, the function evaluations are computed with fixed $M$ and $\alpha$ values rather than Bayesian optimization.  As summarized in Figure~\ref{fig:var}, our robust optimization approach can yield consistently improved lift-to-drag ratios while maintaining a low standard deviation across the tested cases. This demonstrates the effectiveness of the proposed method in achieving aerodynamic performance resilient to variations in operating conditions.

As for the efficiency of the proposed methodology, we consider comparison to a standard technique.  For this purpose, similar to the experiment presented in Figure~\ref{fig:var}, we conduct a multipoint run and solve \eqref{prob:gro} with $c_R(x)=c(x;p)$ for several realizations of $p\in U_p$.  However, this time we solve the resulting problems via the derivative-based solver SLSQP by providing adjoint evaluations.  This is possible as each realization of $p$ produces a deterministic problem for which the adjoint information can be computed via SU2.  The output of this experiment is given in Appendix A.  We observe that the total number of evaluations (objective function and adjoint evaluations) required by each run of SLSQP is in the interval $[32,102]$.  Recall that the optimization process with GLORO has required a total of 217 evaluations, which together with the reported time values imply that within the computation bugdet required by GLORO the multipoint approach can complete at most 3-4 solves.  We anticipate that the efficiency advantage of GLORO would be more prominent when the number of uncertain parameters (i.e. the size of $p$) is larger than two, since the number of grid points increases exponentially for a multipoint solve.
 
We also attempted to gather the results in the literature reported for relevant optimization problems.  Table~\ref{tab:literature_comparison} in Appendix B presents a summary.  There is no much information reported on the computational cost required by the different techniques suggested in these work.  From the few examples, we can deduce that the cost of global-model-based approaches can be relatively high.

\bigskip

\section{Conclusions and Future Work}\label{sec:future}

In this paper we proposed a new methodology that we named GLORO for robust optimization with expensive function evaluations.  The new methodology uses global surrogates for uncertainty quantification (lower level), which is rebuilt for each design obtained during the optimization process.  These surrogates in the parameter space are constructed during the process of Bayesian optimization, which provides sample efficiency.  The upper level optimization occurs in the variable space, and is based on a local model-based DFO technique.  With the motivation of requiring few function evaluations, our model-based DFO algorithm is based on random sampling with a consideration of well-poisedness. 

The design of the methodology was motivated by the robust airfoil design optimization problem, and it was applied to optimize the RAE2822 airfoil under operational uncertainties, with the Mach number varying within [0.715, 0.735] and the angle of attack within [2.82, 3.02] degrees. The objective was to maximize the lift-to-drag ratio across this operational domain. The results confirm the effectiveness of our robust approach, which yielded consistently improved performance in low-performing cases, including a maximum increase of 52\% in the lift-to-drag ratio for specific flight conditions.  Furthermore, a low standard deviation for this ratio is maintained across the entire uncertainty range, demonstrating the robustness of the optimized design.  For this problem with 30 design variables, the optimization process was completed within only a total of 217 CFD simulations. 

The approximate robust function evaluations returned by the Bayesian optimization, and the randomness included in the steps of the inner and outer optimization techniques raise up the question of sensitivity: How large the results provided by GLORO are effected by the randomness and noise?  This is a question to be addresses in our future research.  %Another interesting topic to study is the complexity of models constructed in outer iterations of GLORO: Could it work better if we constructed higher order models in potentially larger trust regions?   
Another issue regarding sensitivity is the choice of the algorithm parameters.  The sensitivity of the results to the parameters of GLORO as well as the parameters of Bayesian Optimization algorithm can be observed through a comprehensive experimental study.  The mathematical analysis of GLORO as an optimization method working in the presence of stochastic noise is an intriguing topic we are planning to work on.  

Another direction for future research is to test GLORO on different airfoil design optimization problems.  We believe a compelling extension of this work would be to apply our robust optimization methodology to a different airfoil geometry and a distinct flight regime, defined by a new Mach number and angle of attack, with new uncertainty ranges established for these parameters. Such a study would powerfully demonstrate the broad applicability of the presented approach. %We believe that an interesting extension of the problem is to consider robustness for variations in design variables as well as the variations in parameters, which could model the uncertainties in manufacturing.

\bigskip

\section*{Data Availability}
The data and the codes used this study are available from the corresponding author upon reasonable request.

%\newpage

\section*{Acknowledgements}
This work is supported by TUSA\c{S} under the LIFT-UP programme with project number 66ab440b93bed. Computing resources used in this work were provided by the National Center for High Performance Computing of Turkey (UHeM) under grant number 4021782025.  Figen Oztoprak and Dilara Bük was supported by TUBITAK under grant number 124M396.  Alpaslan, Dilara and Yusuf were supported by TUBITAK 2209B programme with project number 1139B412403023.

%\newpage
\bibliographystyle{plain}
\bibliography{references}

@article{augustin2017trust,
  title={A trust-region method for derivative-free nonlinear constrained stochastic optimization},
  author={Augustin, F and Marzouk, YM},
  journal={arXiv preprint arXiv:1703.04156},
  year={2017}
}

@article{cozad2014learning,
  title={Learning surrogate models for simulation-based optimization},
  author={Cozad, Alison and Sahinidis, Nikolaos V and Miller, David C},
  journal={AIChE Journal},
  volume={60},
  number={6},
  pages={2211--2227},
  year={2014},
  publisher={Wiley Online Library}
}

@book{conn2009,
  title={Introduction to derivative-free optimization},
  author={Conn, Andrew R and Scheinberg, Katya and Vicente, Luis N},
  year={2009},
  publisher={SIAM}
}

@article{conn2012bilevel,
  title={Bilevel derivative-free optimization and its application to robust optimization},
  author={Conn, Andrew R and Vicente, Lu{\'\i}s Nunes},
  journal={Optimization Methods and Software},
  volume={27},
  number={3},
  pages={561--577},
  year={2012},
  publisher={Taylor \& Francis}
}

@article{catalani2023comparative,
  title={A comparative study of learning techniques for the compressible aerodynamics over a transonic RAE2822 airfoil},
  author={Catalani, Giovanni and Costero, Daniel and Bauerheim, Michael and Zampieri, Luca and Chapin, Vincent and Gourdain, Nicolas and Baqu{\'e}, Pierre},
  journal={Computers \& Fluids},
  volume={251},
  pages={105759},
  year={2023},
  publisher={Elsevier}
}

@article{keane2020robust,
  title={Robust design optimization using surrogate models},
  author={Keane, Andy J and Voutchkov, Ivan I},
  journal={Journal of Computational Design and Engineering},
  volume={7},
  number={1},
  pages={44--55},
  year={2020},
  publisher={Oxford University Press}
}

@article{mulvey1995robust,
  title={Robust optimization of large-scale systems},
  author={Mulvey, John M and Vanderbei, Robert J and Zenios, Stavros A},
  journal={Operations research},
  volume={43},
  number={2},
  pages={264--281},
  year={1995},
  publisher={INFORMS}
}

@article{ben2002robust,
  title={Robust optimization--methodology and applications},
  author={Ben-Tal, Aharon and Nemirovski, Arkadi},
  journal={Mathematical programming},
  volume={92},
  pages={453--480},
  year={2002},
  publisher={Springer}
}

@article{zhang2007general,
  title={General robust-optimization formulation for nonlinear programming},
  author={Zhang, Yin},
  journal={Journal of optimization theory and applications},
  volume={132},
  number={1},
  pages={111--124},
  year={2007},
  publisher={Springer}
}

@book{forrester2008engineering,
  title={Engineering design via surrogate modelling: a practical guide},
  author={Forrester, Alexander and Sobester, Andras and Keane, Andy},
  year={2008},
  publisher={John Wiley \& Sons}
}

@article{sabater2022robust,
  title={Robust design of transonic natural laminar flow wings under environmental and operational uncertainties},
  author={Sabater, Christian and Bekemeyer, Philipp and G{\"o}rtz, Stefan},
  journal={AIAA Journal},
  volume={60},
  number={2},
  pages={767--782},
  year={2022},
  publisher={American Institute of Aeronautics and Astronautics}
}

@inproceedings{sudret2017surrogate,
  title={Surrogate models for uncertainty quantification: An overview},
  author={Sudret, Bruno and Marelli, Stefano and Wiart, Joe},
  booktitle={2017 11th European conference on antennas and propagation (EUCAP)},
  pages={793--797},
  year={2017},
  organization={IEEE}
}

@article{dellino2015metamodel,
  title={Metamodel-based robust simulation-optimization: An overview},
  author={Dellino, Gabriella and Kleijnen, Jack PC and Meloni, Carlo},
  journal={Uncertainty Management in Simulation-Optimization of Complex Systems: Algorithms and Applications},
  pages={27--54},
  year={2015},
  publisher={Springer}
}

@article{bertsimas2010nonconvex,
  title={Nonconvex robust optimization for problems with constraints},
  author={Bertsimas, Dimitris and Nohadani, Omid and Teo, Kwong Meng},
  journal={INFORMS journal on computing},
  volume={22},
  number={1},
  pages={44--58},
  year={2010},
  publisher={INFORMS}
}

@article{gorissen2015practical,
  title={A practical guide to robust optimization},
  author={Gorissen, Bram L and Yanikoglu, Ihsan and Den Hertog, Dick},
  journal={Omega},
  volume={53},
  pages={124--137},
  year={2015},
  publisher={Elsevier}
}

@article{kolvenbach2018approach,
  title={An approach for robust PDE-constrained optimization with application to shape optimization of electrical engines and of dynamic elastic structures under uncertainty},
  author={Kolvenbach, Philip and Lass, Oliver and Ulbrich, Stefan},
  journal={Optimization and Engineering},
  volume={19},
  pages={697--731},
  year={2018},
  publisher={Springer}
}

@article{papadimitriou2016aerodynamic,
  title={Aerodynamic shape optimization for minimum robust drag and lift reliability constraint},
  author={Papadimitriou, Dimitrios I and Papadimitriou, Costas},
  journal={Aerospace Science and Technology},
  volume={55},
  pages={24--33},
  year={2016},
  publisher={Elsevier}
}

@article{cook2017robust,
  title={Robust airfoil optimization and the importance of appropriately representing uncertainty},
  author={Cook, Laurence W and Jarrett, Jerome P},
  journal={AIAA Journal},
  volume={55},
  number={11},
  pages={3925--3939},
  year={2017},
  publisher={American Institute of Aeronautics and Astronautics}
}

@article{ng2014multifidelity,
  title={Multifidelity approaches for optimization under uncertainty},
  author={Ng, Leo WT and Willcox, Karen E},
  journal={International Journal for numerical methods in Engineering},
  volume={100},
  number={10},
  pages={746--772},
  year={2014},
  publisher={Wiley Online Library}
}

@article{garreis2017constrained,
  title={Constrained optimization with low-rank tensors and applications to parametric problems with PDEs},
  author={Garreis, Sebastian and Ulbrich, Michael},
  journal={SIAM Journal on Scientific Computing},
  volume={39},
  number={1},
  pages={A25--A54},
  year={2017},
  publisher={SIAM}	
}

@misc{SU2rae2822,
 title={Resources for the tutorial Turbulent\_2D\_Constrained\_RAE2822},
note={\texttt{https://github.com/su2code/Tutorials/tree/master/design/ Turbulent\_2D\_Constrained\_RAE2822}}
}

@article{su2,
  title={SU2: An open-source suite for multiphysics simulation and design},
  author={Economon, Thomas D and Palacios, Francisco and Copeland, Sean R and Lukaczyk, Trent W and Alonso, Juan J},
  journal={Aiaa Journal},
  volume={54},
  number={3},
  pages={828--846},
  year={2016},
  publisher={American Institute of Aeronautics and Astronautics}
}

@article{chen2015randomized,
  title={Randomized derivative-free optimization of noisy convex functions},
  author={Chen, Ruobing and Wild, Stefan},
  journal={arXiv preprint arXiv:1507.03332},
  year={2015}
}

@article{chen2018stochastic,
  title={Stochastic optimization using a trust-region method and random models},
  author={Chen, Ruobing and Menickelly, Matt and Scheinberg, Katya},
  journal={Mathematical Programming},
  volume={169},
  number={2},
  pages={447--487},
  year={2018},
  publisher={Springer}
}

@article{nesterov2017random,
  title={Random gradient-free minimization of convex functions},
  author={Nesterov, Yurii and Spokoiny, Vladimir},
  journal={Foundations of Computational Mathematics},
  volume={17},
  number={2},
  pages={527--566},
  year={2017},
  publisher={Springer}
}

@article{menickelly2023avoiding,
  title={Avoiding geometry improvement in derivative-free model-based methods via randomization},
  author={Menickelly, Matt},
  journal={arXiv preprint arXiv:2305.17336},
  year={2023}
}

@article{hare2025expected,
  title={Expected decrease for derivative-free algorithms using random subspaces},
  author={Hare, Warren and Roberts, Lindon and Royer, Cl{\'e}ment},
  journal={Mathematics of Computation},
  volume={94},
  number={351},
  pages={277--304},
  year={2025}
}

@article{cook1979aerofoil,
  title={Aerofoil rae 2822-pressure distributions, and boundary layer and wake measurements. experimental data base for computer program assessment},
  author={Cook, PH and McDonald, MA and Firmin, MCP},
  journal={AGARD report ar},
  volume={138},
  pages={47},
  year={1979}
}

@article{garbaruk2003numerical,
  title={Numerical study of wind-tunnel walls effects on transonic airfoil flow},
  author={Garbaruk, Andrey and Shur, Mikhail and Strelets, Mikhail and Spalart, Philippe R},
  journal={AIAA journal},
  volume={41},
  number={6},
  pages={1046--1054},
  year={2003}
}

@article{da2020sensitivity,
  title={Sensitivity and calibration of turbulence model in the presence of epistemic uncertainties},
  author={Da Ronch, Andrea and Panzeri, Marco and Drofelnik, Jernej and d’Ippolito, Roberto},
  journal={CEAS Aeronautical Journal},
  volume={11},
  number={1},
  pages={33--47},
  year={2020},
  publisher={Springer}
}

@article{zhou1995turbulent,
  title={Turbulent transonic airfoil flow simulation using a pressure-based algorithm},
  author={Zhou, Gang and Davidson, Lars and Olsson, Erik},
  journal={AIAA Journal},
  volume={33},
  number={1},
  pages={42--47},
  year={1995}
}

@article{chen2024wind,
  title={Wind tunnel wall interference correction for transonic airfoils with data-reduced ensemble Kalman filter},
  author={Chen, Xin and Wang, Gang and Ye, Zhengyin},
  journal={Physics of Fluids},
  volume={36},
  number={10},
  year={2024},
  publisher={AIP Publishing}
}

@article{zhang2020inverse,
  title={An inverse design method for airfoils based on pressure gradient distribution},
  author={Zhang, Yufei and Yan, Chongyang and Chen, Haixin},
  journal={Energies},
  volume={13},
  number={13},
  pages={3400},
  year={2020},
  publisher={MDPI}
}

@article{chang1988further,
  title={Further comparisons of interactive boundary-layer and thin-layer Navier-Stokes procedures},
  author={Chang, Kuei-Chung and Alemdaroglu, N and Mehta, Unmeel and Cebeci, Tuncer},
  journal={Journal of aircraft},
  volume={25},
  number={10},
  pages={897--903},
  year={1988}
}

@inproceedings{weinmeister2024,
  title={Comparison of Full-Field and Integrated CFD Convergence Based on Richardson Extrapolation},
  author={Weinmeister, Justin R and Sanjaya, Devina P},
  booktitle={AIAA AVIATION FORUM AND ASCEND 2024},
  pages={4470},
  year={2024}
}

@article{vuruskan2019,
  title={Impact of turbulence models and shape parameterization on robust aerodynamic shape optimization},
  author={Vuruskan, Aslihan and Hosder, Serhat},
  journal={Journal of Aircraft},
  volume={56},
  number={3},
  pages={1099--1115},
  year={2019},
  publisher={American Institute of Aeronautics and Astronautics}
}

@article{jinks2018optimisation,
  title={Optimisation of adaptive shock control bumps with structural constraints},
  author={Jinks, Edward and Bruce, Paul and Santer, Matthew},
  journal={Aerospace Science and Technology},
  volume={77},
  pages={332--343},
  year={2018},
  publisher={Elsevier}
}

@article{sabater2020efficient,
  title={Efficient bilevel surrogate approach for optimization under uncertainty of shock control bumps},
  author={Sabater, Christian and Bekemeyer, Philipp and G{\"o}rtz, Stefan},
  journal={AIAA Journal},
  volume={58},
  number={12},
  pages={5228--5242},
  year={2020},
  publisher={American Institute of Aeronautics and Astronautics}
}

@article{mazaheri2015,
  title={Optimization and analysis of shock wave/boundary layer interaction for drag reduction by shock control bump},
  author={Mazaheri, K and Kiani, KC and Nejati, A and Zeinalpour, M and Taheri, R},
  journal={Aerospace Science and Technology},
  volume={42},
  pages={196--208},
  year={2015},
  publisher={Elsevier}
}

@article{hock1980,
  title={Test examples for nonlinear programming codes},
  author={Hock, Willi and Schittkowski, Klaus},
  journal={Journal of optimization theory and applications},
  volume={30},
  number={1},
  pages={127--129},
  year={1980},
  publisher={Springer}
}

@article{nemati2020robust,
  title={Robust aerodynamic morphing shape optimization for high-lift missions},
  author={Nemati, M and Jahangirian, A},
  journal={Aerospace Science and Technology},
  volume={103},
  pages={105897},
  year={2020},
  publisher={Elsevier}
}

@article{tao2019application,
  title={Application of deep learning based multi-fidelity surrogate model to robust aerodynamic design optimization},
  author={Tao, Jun and Sun, Gang},
  journal={Aerospace Science and Technology},
  volume={92},
  pages={722--737},
  year={2019},
  publisher={Elsevier}
}

@article{fusi2020assessment,
  title={Assessment of robust optimization for design of rotorcraft airfoils in forward flight},
  author={Fusi, Francesca and Congedo, Pietro Marco and Guardone, Alberto and Quaranta, Giuseppe},
  journal={Aerospace Science and Technology},
  volume={107},
  pages={106355},
  year={2020},
  publisher={Elsevier}
}

@inproceedings{jameson1981,
  title={Numerical solution of the Euler equations by finite volume methods using Runge Kutta time stepping schemes},
  author={Jameson, Antony and Schmidt, Wolfgang and Turkel, Eli},
  booktitle={14th fluid and plasma dynamics conference},
  pages={1259},
  year={1981}
}

@book{garnett2023,
  title={Bayesian optimization},
  author={Garnett, Roman},
  year={2023},
  publisher={Cambridge University Press}
}

@article{bandeira2014,
  title={Convergence of trust-region methods based on probabilistic models},
  author={Bandeira, Afonso S and Scheinberg, Katya and Vicente, Luis Nunes},
  journal={SIAM Journal on Optimization},
  volume={24},
  number={3},
  pages={1238--1264},
  year={2014},
  publisher={SIAM}
}

@article{maggiar2018,
  title={A derivative-free trust-region algorithm for the optimization of functions smoothed via gaussian convolution using adaptive multiple importance sampling},
  author={Maggiar, Alvaro and Wachter, Andreas and Dolinskaya, Irina S and Staum, Jeremy},
  journal={SIAM Journal on Optimization},
  volume={28},
  number={2},
  pages={1478--1507},
  year={2018},
  publisher={SIAM}
}

@article{croicu2012robust,
  title={Robust airfoil optimization using maximum expected value and expected maximum value approaches},
  author={Croicu, Ana-Maria and Hussaini, M Yousuff and Jameson, Antony and Klopfer, Goetz},
  journal={AIAA journal},
  volume={50},
  number={9},
  pages={1905--1919},
  year={2012}
}

@article{dodson2009robust,
  title={Robust aerodynamic design optimization using polynomial chaos},
  author={Dodson, Michael and Parks, Geoffrey T},
  journal={Journal of Aircraft},
  volume={46},
  number={2},
  pages={635--646},
  year={2009}
}

@article{ji2025d,
  title={D-optimal polynomial chaos expansion for adjoint-based aerodynamic robust optimization in transonic flows},
  author={Ji, Xinze and Yang, Tihao and Shi, Yayun and Ma, Yuhang and Bai, Junqiang},
  journal={Aerospace Science and Technology},
  pages={110659},
  year={2025},
  publisher={Elsevier}
}

@article{keane2012cokriging,
  title={Cokriging for robust design optimization},
  author={Keane, Andy J},
  journal={AIAA journal},
  volume={50},
  number={11},
  pages={2351--2364},
  year={2012}
}

@article{zhang2024efficient,
  title={An efficient robust aerodynamic design optimization method based on a multi-level hierarchical Kriging model and multi-fidelity expected improvement},
  author={Zhang, Yu and Han, Zhong-hua and Song, Wen-ping},
  journal={Aerospace Science and Technology},
  volume={152},
  pages={109401},
  year={2024},
  publisher={Elsevier}
}

\newpage
\section*{Appendix A.  Computational Cost of A Multipoint Run}

The cost of a multipoint run with SLSQP employing adjoint evaluations is summarized in Table~\ref{tab:costmultipoint} as a reference to assess the efficiency of GLORO.  Each row reports the computational cost of an individual solve within the multipoint run.  %Each solution obtained with a fixed $p$ value is evaluated across all 16 realizations within the uncertainty box; maximum, minimum, mean and standard deviation statistics are computed on the resulting 16 $c_L/c_D$ values.

\bigskip

\begin{table}[H]
\caption{Cost of multipoint SLSQP solve}
\label{tab:costmultipoint}
%\footnotesize
%\begin{tabular}{r | c c c c c |c c c c c }
%\toprule
%Fixed value &	&	SLSQP	&	$\#$ of	&	$\#$	of &		Final &	\multicolumn{4}{c}{$c_L/c_D$ for other $p$}	\\
%of $p = (M,a)$	&	Exit	&	Iters	&	Primal	&	Adjoint	&	 $c_L/c_D$	&	Mean	&	Min	&	Max	&	St.Dev.	\\
%\midrule
%(0.715,	2.820)	&	0	&	18	&	48	&	18	&	81.28	&	60.52	&	43.46	&	81.28	&	12.16	\\
%(0.715,	2.887)	&	0	&	7	&	25	&	7	&	76.20	&	59.42	&	43.54	&	76.17	&	10.55	\\
%(0.715,	2.953)	&	0	&	27	&	50	&	27	&	82.97	&	63.86	&	47.35	&	82.97	&	10.96	\\
%(0.715,	3.020)	&	0	&	30	&	57	&	30	&	82.73	&	63.12	&	48.41	&	82.73	&	9.6	\\
%(0.722,	2.820)	&	0	&	35	&	67	&	35	&	80.71	&		&		&		&		\\
%(0.722,	2.887)	&	0	&	34	&	62	&	34	&	80.6	&	65.45	&	53.36	&	79.54	&	6.99	\\
%(0.722,	2.953)	&	0	&	28	&	59	&	28	&	79.85	&	65.05	&	55.01	&	78.78	&	6.09	\\
%(0.722,	3.020)	&	0	&	29	&	48	&	29	&	80.16	&	63.99	&	57.53	&	72.89	&	4.84	\\
%(0.728,	2.820)	&	0	&	32	&	53	&	32	&	78.4	&	67.13	&	59.28	&	78.3	&	5.38	\\
%(0.728,	2.887)	&	0	&	35	&	66	&	35	&	78.16	&	65.55	&	58.77	&	77.97	&	5.53	\\
%(0.728,	2.953)	&	0	&	27	&	62	&	27	&	77.3	&	64.5	&	56.68	&	77.18	&	5.47	\\
%(0.728,	3.020)	&	0	&	33	&	65	&	33	&	77.73	&	62.84	&	53.44	&	77.5	&	5.56	\\
%(0.735,	2.820)	&	0	&	30	&	64	&	30	&	75.48	&	64.73	&	59.03	&	75.48	&	4.96	\\
%(0.735,	2.887)	&	0	&	24	&	53	&	24	&	74.58	&	62.68	&	55.9	&	74.52	&	5.19	\\
%(0.735,	2.953)	&	0	&	9	&	30	&	9	&	68.72	&	52.63	&	45.7	&	68.74	&	6.02	\\
%(0.735,	3.020)	&	0	&	7	&	26	&	7	&	64.63	&	53.31	&	47.9	&	64.63	&	3.6	\\
%\midrule
%		&		&	TOTAL	&	709	&	342	&		&		&		&		&		\\
%\bottomrule
%\end{tabular}
\small
\centering
\begin{tabular}{r | c c c c | c  }
\toprule
Fixed value &	&	$\#$ of	&	$\#$ of	&	$\#$	of &	CPU	 \\ 
of $p = (M,a)$	&	Exit	&	Iterations	&	Primal Eval.	&	Adjoint Eval.	&  hours\\
\midrule
(0.715,	2.820)	&	0	&	18	&	48	&	18	&	13.94\\
(0.715,	2.887)	&	0	&	7	&	25	&	7	&	6.16\\
(0.715,	2.953)	&	0	&	27	&	50	&	27	&	18.28\\
(0.715,	3.020)	&	0	&	30	&	57	&	30	&	20.35\\
(0.722,	2.820)	&	0	&	35	&	67	&	35	&	23.23 \\
(0.722,	2.887)	&	0	&	34	&	62	&	34	&	22.80\\
(0.722,	2.953)	&	0	&	28	&	59	&	28	&	19.86\\
(0.722,	3.020)	&	0	&	29	&	48	&	29	&	18.85\\
(0.728,	2.820)	&	0	&	32	&	53	&	32	&	20.84\\
(0.728,	2.887)	&	0	&	35	&	66	&	35	&      23.71\\
(0.728,	2.953)	&	0	&	27	&	62	&	27	&      19.75\\
(0.728,	3.020)	&	0	&	33	&	65	&	33	&      22.77\\
(0.735,	2.820)	&	0	&	30	&	64	&	30	&      21.29\\
(0.735,	2.887)	&	0	&	24	&	53	&	24	&      17.22\\
(0.735,	2.953)	&	0	&	9	&	30	&	9	&	7.68\\
(0.735,	3.020)	&	0	&	7	&	26	&	7	&	6.30\\
\midrule
\rowcolor{gray!12}
GLORO & - & 10 & 217 & 0 & 26.86\\
\bottomrule
\end{tabular}
\end{table}
\clearpage

\newgeometry{a4paper,left=1in,right=1in,top=1in,bottom=1in,nohead}
%\begin{landscape}
\section*{Appendix B.  Summary Information on Relevant Work} %in the Literature
%\begin{table}[H]
\begin{sidewaystable}[h]
\centering
\captionsetup{width=0.95\linewidth}
\caption{Problem formulations, solution methods, and computational
effort in robust aerodynamic shape optimization studies
(``n.r.'' indicates that the required information is not reported).}
\label{tab:literature_comparison}

\scriptsize
\setlength{\tabcolsep}{3pt}
\renewcommand{\arraystretch}{1.08}

\resizebox{0.95\linewidth}{!}{%
\begin{tabular}{@{}
>{\raggedright\arraybackslash}p{2.2cm}
>{\raggedright\arraybackslash}p{4.0cm}
>{\raggedright\arraybackslash}p{3.1cm}
>{\raggedright\arraybackslash}p{3.1cm}
>{\centering\arraybackslash}p{1.3cm}
>{\raggedright\arraybackslash}p{5.3cm}
>{\raggedright\arraybackslash}p{5.5cm}
@{}}
\toprule

\textbf{Study} &
\textbf{Objective} &
\textbf{Uncertainties} &
\textbf{Constraints} &
\textbf{$n$} &
\textbf{Method} &
\textbf{Evaluations} \\

\midrule

Dodson \& Parks (2009) &
Bi-objective: max.\ mean and min.\ variance of $C_L/C_D$ &
Leading-edge thickness &
n.r. &
2 &
Polynomial-chaos UQ\textsuperscript{a} with differential evolution
for RO\textsuperscript{b} &
3{,}000 nonrobust;
8{,}000 / 12{,}000 / 32{,}000 / 16{,}000 robust XFOIL evaluations \\

\addlinespace[3pt]

Croicu et al.\ (2012) &
Single-objective: expected optimal shape (EMV) or mean $C_D$ (MEV) &
Mach &
$C_L$, thickness &
65 &
Monte Carlo EMV/MEV optimization with adjoint gradients;
SYN83 steady inviscid Euler &
Total n.r.; 20 design cycles and 1{,}000 Mach samples
per stochastic method \\

\addlinespace[3pt]

Keane (2012) &
Bi-objective: min.\ mean and std.\ dev.\ of normalized pressure loss &
Manufacturing variation, foreign-object damage, and flank erosion
(16 noise variables) &
Variable bounds &
10 &
Direct and Kriging/co-Kriging-assisted
NSGA-II\textsuperscript{c} with Monte Carlo UQ &
11{,}250 direct; 4{,}875 Kriging; 4{,}875 co-Kriging
HYDRA RANS CFD runs per search \\

\addlinespace[3pt]

Fusi et al.\ (2015) &
Bi-objective: max.\ mean and min.\ variance of $C_L/C_D$;
hovering rotor airfoil &
Blade pitch angle, induced velocity &
Camber-line geometry; pitching moment limit &
17 &
MF\textsuperscript{d} RO using NSGA-II and polynomial-chaos UQ &
Full HF\textsuperscript{e}: 32{,}450--33{,}475;
MF: 6{,}325--10{,}324 HF $+$
29{,}166--31{,}625 LF\textsuperscript{g}
(HF: Euler/IBL\textsuperscript{f}; LF: panel/IBL)\\

\addlinespace[3pt]

Papadimitriou \& Papadimitriou (2016) &
Single-objective weighted formulation: mean $+$ std.\ dev.\ of $C_D$ &
Mach, AoA and geometry &
Lift reliability chance constraint &
23\textsuperscript{h} &
Adjoint-based RO with sparse-grid UQ and
FORM\textsuperscript{i} &
Total n.r.; Euler flow and Euler adjoint evaluations used for
UQ and reliability analysis \\

\addlinespace[3pt]

Cook \& Jarrett (2017) &
Bi-objective: mean and std.\ dev.\ of $C_D$ &
Mach &
$C_L$, geometry, variable bounds &
14 &
Polynomial-chaos UQ with NSGA-II; SU2 RANS &
Total n.r.\ (RANS CFD evaluations) \\

\addlinespace[3pt]

Tao \& Sun (2019) &
Single-objective weighted formulation: squared mean $+$ variance
of $C_D$ &
Mach &
AoA limit; airfoil/section thickness &
12 &
Deep-learning-based MF surrogate optimization using improved PSO &
4{,}250 HF RANS CFD computations and
1{,}600 MF surrogate computations \\

\addlinespace[3pt]

Vuruskan \& Hosder (2019) &
Single-objective weighted formulation: normalized
$\mu_{C_D}$, $\sigma_{C_D}^{2}$, and squared deviation from
target $\mu_{C_L}$ &
Mach &
Area and variable bounds &
n.r.\textsuperscript{j} &
Adjoint-based RO using point-collocation NIPC and SLSQP;
SU2 RANS with SA\textsuperscript{k} and SST\textsuperscript{l} &
Total n.r.; per optimization iteration:
6 direct RANS flow evaluations and 6 RANS adjoint evaluations \\

\addlinespace[3pt]

Sabater et al.\ (2020) &
Single-objective cases: 50\,\% or 95\,\% quantile of $C_D$ &
Mach, $C_L$, bump height &
Constant lift, bump--flap clearance, variable bounds &
5 &
Bilevel GP-assisted RO &
3{,}219 RANS CFD evaluations \\

\addlinespace[3pt]

Keane \& Voutchkov (2020) &
Bi-objective: mean and std.\ dev.\ of normalized pressure loss &
Manufacturing variation, foreign-object damage, and flank erosion
(16 noise variables) &
Unconstrained &
10 &
Co-Kriging, combined Kriging, and combined co-Kriging-assisted
NSGA-II &
3{,}140 co-Kriging; 3{,}150 combined Kriging;
3{,}130 combined co-Kriging HYDRA RANS CFD runs;
about 47{,}000 for the reference Pareto-front estimate \\

\addlinespace[3pt]

Sabater et al.\ (2022) &
Single-objective cases: mean or 95\,\% quantile of $C_D$ &
Critical TS/CF $N$-factors; additionally Mach and $C_L$ &
Constant lift, minimum thickness, variable bounds &
10 &
Bilevel GP-assisted RO with transition prediction &
462 RANS/transition evaluations for deterministic optimization;
200 DoE samples; RO total n.r. \\

\addlinespace[3pt]

Zhang et al.\ (2024) &
Single-objective: mean $+$ std.\ dev.\ of $C_D$ &
Mach, AoA &
Mean $C_L$, mean $C_M$, area, and variable bounds &
18 &
MF RO using MHK\textsuperscript{m} and
MFEI\textsuperscript{n}; RANS with SA &
HF/MiF/LF RANS evaluations:
MFRDO, 180/870/2{,}450 (317.5 eq.\ HF);
progressive, 4/606/2{,}010 (93.6 eq.\ HF)\\

\addlinespace[3pt]
\midrule

\rowcolor{gray!12}
GLORO (present) &
Single-objective worst-case optimization of $C_L/C_D$ &
Mach, AoA &
Min.\ thickness; bounds; $C_M$ monitored &
38 &
Bayesian optimization (GP) in parameter space $+$ local
model-based DFO in design space &
217 RANS primal evaluations and 0 adjoint evaluations;
10 outer iterations \\

\bottomrule
\end{tabular}%
}

\vspace{0.35em}

\begin{minipage}{0.95\linewidth}
\scriptsize
$n$: number of design variables.

\textsuperscript{a}~UQ: uncertainty quantification;\quad
\textsuperscript{b}~RO: robust optimization;\quad
\textsuperscript{c}~NSGA: nondominated sorting genetic algorithm;\quad
\textsuperscript{d}~MF: multifidelity;\quad
\textsuperscript{e}~HF: high fidelity;\quad
\textsuperscript{f}~IBL: integral boundary layer;\quad
\textsuperscript{g}~LF: low fidelity;\quad
\textsuperscript{h}~22 movable B\'ezier control-point coordinates
plus mean angle of attack;\quad
\textsuperscript{i}~FORM: first-order reliability method;\quad
\textsuperscript{j}~For RO, the paper reports 20 B-spline control
points, 256 Hicks--Henne bump-function amplitudes, and a
$2{\times}20$ FFD control lattice, respectively; angle of attack is
an additional design variable; combined totals are not explicitly
reported;\quad
\textsuperscript{k}~SA: Spalart--Allmaras;\quad
\textsuperscript{l}~SST: shear-stress transport;\quad
\textsuperscript{m}~MHK: multi-level hierarchical Kriging;\quad
\textsuperscript{n}~MFEI: multi-fidelity expected improvement;\quad
DoE: design of experiments;\quad
NIPC: non-intrusive polynomial chaos;\quad
SLSQP: sequential least-squares programming;\quad
GP: Gaussian process;\quad
DFO: derivative-free optimization.
\end{minipage}
\end{sidewaystable}
%\end{table}
%\end{landscape}
\restoregeometry

\end{document}